\documentclass[12pt]{amsart}
\usepackage{amsfonts,amssymb,amsthm,eucal}

\newtheorem{thm}{Theorem}

\newtheorem{lemma}[thm]{Lemma}
\newtheorem{prop}[thm]{Proposition}
\newtheorem{claim}{Claim}

\newcommand{\R}{\mathbb{R}}
\newcommand{\E}{\mathbb{E}}
\newcommand{\N}{\mathbb{N}}
\newcommand{\K}{\mathcal{K}}
\DeclareMathOperator{\sgn}{sgn}
\DeclareMathOperator{\vol}{vol}
\DeclareMathOperator{\conv}{conv}

\title{Volumes of symmetric random polytopes}

\author[M. Meckes]{Mark W. Meckes}

\address{Department of Mathematics, Case Western Reserve University,
    Cleveland, Ohio 44106.}

\email{mwm2@po.cwru.edu}

\thanks{Partially supported by a grant from the National Science Foundation.}


\allowdisplaybreaks[2]

\begin{document}

\begin{abstract}
We consider the moments of the volume of the symmetric convex hull of
independent random points in an $n$-dimensional symmetric convex body. We
calculate explicitly the second and fourth moments for $n$ points when the
given body is $B_q^n$ (and all of the moments for the case $q = 2$), and
derive from these the asymptotic behavior of the expected volume of a
random simplex in those bodies.
\end{abstract}

\maketitle

\section{Introduction} \label{S:intro}

Let $\K^n$ denote the family of all convex bodies in the Euclidean space
$\R^n$, that is, all compact convex sets with interior points; and let
$\K_s^n$ denote the family of symmetric convex bodies, that is, all $K \in
\K^n$ such that $K = -K$. For any $K \in \K^n$ and $N \ge n + 1$, we
define the random variable
$$U_{K, N} = \frac{1}{|K|} | \conv \{ x_1, \ldots, x_N \} |,$$
where $|A|$ denotes the volume of a Borel set $A \subset \R^n$, and $x_1,
\ldots, x_N$ are independent random points uniformly distributed in $K$.
That is, $U_{K, N}$ is the normalized volume of a random polytope in $K$;
in particular, $U_K = U_{K, n+1}$ is the normalized volume of a random
simplex in $K$. Note that the distribution of $U_{K, N}$ is an affine
invariant of $K$. For $K \in \K_s^n$ and $N \ge n$, we define
$$V_{K, N} = \frac{1}{|K|} | \conv \{ \pm x_1, \ldots, \pm x_N \} |,$$
where again $x_1, \ldots, x_N$ are independent uniform random points in
$K$. $V_{K, N}$ is the normalized volume of a symmetric random polytope in
$K$; $V_K = V_{K, n}$ is the normalized volume of a random crosspolytope
in $K$. The distribution of $V_{K, N}$ is a linear invariant of $K$.

Interest in the behavior of the moments of $U_{K, N}$ dates back to the
four point problem of Sylvester, first stated in the 1860's (see
\cite{Pfiefer} for a review of this problem's early history). In this
paper we consider also the moments of $V_{K, N}$, and their relationship
to the moments of $U_{K, N}$. We first observe that these moments are
always minimized when $K$ is an ellipsoid.

\begin{thm} \label{T:min}
Let $f : \R_+ \to \R_+$ be a strictly increasing continuous function.
\begin{enumerate}
\item \label{I:nonsymm-min} $\E f(U_{K, N}) \ge \E f(U_{B_2^n, N})$ for
any $K \in \K^n$ and $N \ge n + 1$.
\item \label{I:symm-min} $\E f(V_{K, N}) \ge \E f(V_{B_2^n, N})$ for any
$K \in \K_s^n$ and $N \ge n$.
\end{enumerate}
\end{thm}

Special cases of Theorem \ref{T:min}(\ref{I:nonsymm-min}) were proved in
\cite{Blaschke1,Blaschke2,Groemer1,Groemer2,Schoepf}; the version quoted
here was proved by Giannopoulos and Tsolomitis in \cite{GT}, and was
extended by Hartzoulaki and Paouris to quermassintegrals of random
polytopes in \cite{HP}. The main tool in all of the proofs is Steiner
symmetrization, and only notational changes in the proofs are necessary to
prove Theorem \ref{T:min}(\ref{I:symm-min}).

A natural conjecture is that the moments $\E U_{K, N}^p$, $p \ge 1$, are
always maximized for all $K \in \K^n$ when $K$ is a simplex, and that $\E
U_{K, N}^p$ and $\E V_{K, N}^p$ are maximized for all $K \in \K_s^n$ when
$K$ is a parallelotope or crosspolytope. Both of these conjectures are
known to be true when $n = 2$, but only partial results in this direction
are known for $n \ge 3$; see \cite{DL,Giann,CCG,Meckes3}.

For the rest of this paper, we will consider only the particular cases
$U_K = U_{K, n+1}$ and $V_K = V_{K, n}$, that is, the volume of a random
$n$-dimensional simplex or crosspolytope in $K \in \K_s^n$. In Section
\ref{S:gen}, we make some general observations about the moments of $U_K$
and $V_K$; in particular, we show that $\| V_K \|_p^{1/n} \simeq 2 \|
U_K\|_p^{1/n}$ for all $1 \le p \le \infty$, where here and below $f
\simeq g$ means $f = (1 + o(1)) g$ as $n \to \infty$. In Section
\ref{S:computing}, we calculate $\E V_{B_q^n}^2$ and $\E V_{B_q^n}^4$
explicitly, where $B_q^n = \{ x \in \R^n : \sum_{i=1}^n |x_i|^q \le 1 \}$
for $1 \le q < \infty$, and $B_\infty^n = \{ x \in \R^n : \max |x_i| \le 1
\}$. From this, we are able to derive the exact asymptotic order of $(\E
V_{B_q^n})^{1/n}$ and $(\E U_{B_q^n})^{1/n}$. Finally, in Section
\ref{S:euclidean}, we calculate all of the moments of $V_{B_2^n}$.

\section{General observations} \label{S:gen}

Recall that an Orlicz function is a continuous, increasing, convex
function $\psi : \R_+ \to \R_+$ such that $\psi(0) = 0$ and $\lim_{t \to
\infty} \psi(t) = \infty$; the corresponding Orlicz norm of a random
variable $X$ is
$$\| X \|_\psi = \inf \left\{ \rho > 0 : \E \psi \left(\frac{|X|}{\rho}
    \right) \le 1 \right\}.$$

\begin{prop} \label{T:V_K,U_K}
Let $K \in \K_s^n$. Then for any Orlicz function $\psi$, we have
$$\frac{2}{(n+1)^{1/n}}\|U_K\|_\psi^{1/n} \le \|V_K\|_\psi^{1/n}
     \le 2\|U_K\|_\psi^{1/n}.$$
In particular, for any $1 \le p \le \infty$,
$$\frac{2}{(n+1)^{1/n}}\|U_K\|_p^{1/n} \le \|V_K\|_p^{1/n}
     \le 2\|U_K\|_p^{1/n}.$$
\end{prop}
\begin{proof}
Define the functions
\begin{eqnarray*}
f_0(x_0,x_1,\ldots,x_n)&=&\det(x_1,\ldots,x_n), \\
f_i(x_0,x_1,\ldots,x_n)&=&\det(x_1,\ldots,x_{i-1},x_0,x_{i+1},\ldots,x_n)
    \text{ for } i=1, \ldots, n.
\end{eqnarray*}
To simplify notation, we assume that $|K|=1$. Since $\psi$ is convex and
nondecreasing, for any $\rho > 0$,
\begin{eqnarray*}
\E \psi\left(\frac{n! U_K }{\rho}\right) & = &
    \int_{K^{n+1}} \psi\left(\frac{1}{\rho} \left| \sum_{i=0}^n f_i(-x_0, x_1,
        \ldots, x_n) \right| \right) dx_0 \cdots dx_n \\
& \le & \frac{1}{n+1} \sum_{i=0}^n \int_{K^{n+1}}
        \psi\left(\frac{n+1}{\rho} |f_i(-x_0, x_1, \ldots, x_n)| \right) dx_0
        \cdots dx_n \\
& = & \E \psi\left(\frac{(n+1)!}{2^n \rho} V_K \right),
\end{eqnarray*}
which implies that $\| U_K \|_\psi \le \frac{n+1}{2^n} \| V_K \|_\psi$. On
the other side we have
\begin{eqnarray*}
\E \psi\left( \frac{n! V_K }{2^n \rho}\right)
    & = & \int_{K^{n+1}} \psi\left( \frac{1}{2\rho}\left| \sum_{i=0}^n
        f_i(x_0, x_1,\ldots, x_n) + \sum_{i=0}^n f_i(-x_0, x_1, \ldots,
        x_n) \right| \right) dx_0 \cdots dx_n \\
& \le & \frac{1}{2} \int_{K^{n+1}} \psi\left( \frac{1}{\rho} \left|
        \sum_{i=0}^n f_i(x_0, x_1, \ldots, x_n) \right| \right) dx_0
        \cdots dx_n \\*
& & \qquad + \frac{1}{2} \int_{K^{n+1}} \psi\left( \frac{1}{\rho}
        \left| \sum_{i=0}^n f_i(-x_0, x_1, \ldots, x_n) \right| \right)
        dx_0 \cdots dx_n \\
& = & \E \left( \frac{n!}{\rho} U_K \right),
\end{eqnarray*}
where we have also used the symmetry of $K$. This implies that $\| V_K
\|_\psi \le 2^n \| U_K \|_\psi$. The $L_p$ case for $1 \le p < \infty$
follows by setting $\psi(t) = t^p$; the $L_\infty$ case then follows by
letting $p \to \infty$.
\end{proof}

Intuitively, Proposition \ref{T:V_K,U_K} implies that at the scale of
$n$th roots, for large $n$, the volume of a random crosspolytope in $K$ is
behaves similarly to $2^n$ times the volume of a random simplex in $K$.
The proof above generalizes the proof of Proposition 5.6 in \cite{MP},
which is essentially the $L_1$ case of Proposition \ref{T:V_K,U_K} stated
in different terms.

By a standard application of Borell's lemma (see \cite[Appendix III]{MS}),
one also obtains the following.

\begin{prop}\label{T:norm-equiv}
For $1 \le p < q < \infty$, there are constants $c_{p,q} > 0$ such that
$$\| V_K \|_p^{1/n} \le \| V_K \|_q^{1/n}
    \le c_{p,q} \| V_K \|_p^{1/n}$$
for every $K \in \K_s^n$, and
$$\| U_K \|_p^{1/n} \le \| U_K \|_q^{1/n}
    \le c_{p,q} \| U_K \|_p^{1/n}$$
for every $K \in \K^n$.
\end{prop}

The isotropic constant $L_K$ of a body $K \in \K^n$ can be defined by
$$L_K^{2n} = \det \left[ \frac{1}{|K|^{n+2}} \int_K x_i x_j dx \right].$$
By expanding the determinant expressions, one obtains
\begin{equation} \label{E:V_K-isotropic}
\left(\E V_K^2
\right)^{1/n} = \frac{4}{(n!)^{1/n}} L_K^2
    \simeq \frac{4 e}{n} L_K^2
\end{equation}
for any $K \in \K_s^n$, and
$$\left(\E U_K^2 \right)^{1/n} = \left(\frac{n+1}{n!}\right)^{1/n} L_K^2
    \simeq \frac{e}{n} L_K^2$$
for any $K \in \K^n$ with centroid at the origin. It is known that for
some constant $c > 0$, $L_K \le c n^{1/4} \log n$ for any $K \in \K^n$,
due to Bourgain \cite{Bourgain} when $K \in \K_s^n$ and Paouris
\cite{Paouris} in the general case. Combining this with Proposition
\ref{T:norm-equiv}, we conclude that for $p \ge 1$, there are constants
$c_p > 0$ such that
$$ \| U_K \|_p^{1/n} \le c_p n^{-1/4} \log n $$
for every $K \in \K^n$ and
$$ \| V_K \|_p^{1/n} \le c_p n^{-1/4} \log n $$
for every $K \in \K_s^n$. Recall that the well-known hyperplane conjecture
for symmetric convex bodies is equivalent to the conjecture that $L_K =
O(1)$ for $K \in \K_s^n$ (see \cite{MP}). By the above observations, this
is equivalent to the conjecture that for some $p \ge 1$, $\| U_K
\|_p^{1/n}$ or $\| V_K \|_p^{1/n}$ is $O(n^{-1/2})$. We remark also that
since
$$\left( \E U_{\Delta_n}^2 \right)^{1/n} = \frac{(n!)^{1/n}}{(n + 1)(n + 2)}
    \simeq \frac{1}{e n},$$
where $\Delta_n$ denotes an $n$-dimensional simplex, the hyperplane
conjecture is implied by the conjecture that simplices maximize the
moments (or, by Proposition \ref{T:norm-equiv}, just the expected value)
of $U_K$.

\section{Random polytopes in 1-symmetric bodies}\label{S:computing}

In this section we derive an expression for the fourth moment of $V_K$
when $K$ is {\em 1-symmetric}, that is, when $K$ is invariant under all
reflections in the coordinate hyperplanes and all permutations of the
coordinates of $\R^n$. Combining this with the connection between the
second moment and the isotropic constant $L_K$, and Proposition
\ref{T:V_K,U_K}, we are able to determine the asymptotic behavior of $(\E
V_K)^{1/n}$ and $(\E U_K)^{1/n}$ when $K = B_q^n$. The approach here is
inspired by a suggestion of Kingman \cite{Kingman}, who noted that since
$U_K$ is a bounded nonnegative random variable, its distribution is
uniquely determined by its even moments, so that one might study that
distribution for simple enough bodies by expanding the determinant
expressions for those moments and integrating explicitly. Here we apply
this approach not to $U_K$ directly, but to $V_K$, which is more tractable
because of the simpler determinant expression for the volume of the
symmetric convex hull of $n$ points.

\begin{prop} \label{T:moment4}
Suppose $n \ge 2$ and $K \in \K_s^n$ is 1-symmetric. Then
$$\E V_K^4 = \frac{16^n (n+1)(n+2)}{2 (n!)^2 |K|^{n+4}}
    \left( \int_K x_1^2 x_2^2 dx \right)^n \varphi(A_K - 3,n),$$
where
\begin{equation} \label{E:phi}
\varphi(t,n) = \left( 1 - \frac{2t}{(n+2)} +
    \frac{t^2}{(n+1)(n+2)} \right)
    \sum_{k = 0}^n \frac{t^k}{k!}
    + \frac{t^{n+1}}{(n+1)!} - \frac{t^{n+2}}{(n+2)!}
\end{equation}
and $A_K = \frac{\int_K x_1^4 dx}{\int_K x_1^2 x_2^2 dx}$.
\end{prop}

Note that for $t \in \R$ fixed, $\varphi(t,n) \simeq e^t $. The proof of
Proposition \ref{T:moment4}, which is combinatorial in nature, is
postponed until the end of this section.

By \eqref{E:V_K-isotropic} and Proposition \ref{T:moment4},
\begin{eqnarray*}
\E V_{B_q^n}^2 & = & \frac{1}{n!}
    \left(\frac{\Gamma(1+\frac{3}{q}) \Gamma(1+\frac{n}{q})}
        {3\Gamma(1+\frac{1}{q})^3 \Gamma(1+\frac{n+2}{q})} \right)^n
    \Gamma\left(1+\frac{n}{q}\right)^2, \\
\E V_{B_q^n}^4 & = & \frac{(n+1)(n+2)}{2 (n!)^2}
    \left(\frac{\Gamma(1+\frac{3}{q})^2 \Gamma(1+\frac{n}{q})}
        {9\Gamma(1+\frac{1}{q})^6 \Gamma(1+\frac{n+4}{q})} \right)^n
    \Gamma\left(1+\frac{n}{q}\right)^4 \varphi(A_q - 3,n),
\end{eqnarray*}
where $A_q = A_{B_q^n} = \frac{9 \Gamma(1+1/q) \Gamma(1+5/q)}{5
\Gamma(1+3/q)^2}$. Using Stirling's formula, we have, for all $1 \le q \le
\infty$,
\begin{eqnarray}
\left( \E V_{B_q^n}^2 \right)^{1/n} & \simeq &
    \alpha(1/q) n^{-1}, \label{E:B_q^n-2}\\
\left( \E V_{B_q^n}^4 \right)^{1/n} & \simeq &
    \alpha(1/q)^2 n^{-2}, \label{E:B_q^n-4}
\end{eqnarray}
where $\alpha : [0, 1] \to \R$ is given by
$$\alpha(t) = \frac{e^{1 - 2t} \Gamma(1 + 3t)}{3 \Gamma(1 + t)^3}.$$
Recall that Corollary \ref{T:norm-equiv} implies that $(\E
V_{B_q^n})^{1/n}$ is of the same order as $(\E V_{B_q^n}^2)^{1/2n}$;
however, we can derive a more precise estimate for the expected value by
using the following. (This is an estimate used in the proof of
Khintchine's inequality; see \cite[Theorem 2.b.3]{LT1}.)

\begin{lemma} \label{T:Khintchine}
Let $X$ be a random variable such that $\E |X|^4 < \infty$. Then
$$ \left( \frac{\| X \|_2}{\| X \|_4} \right)^2 \| X \|_2
    \le \| X \|_1 \le \| X \|_2.$$
\end{lemma}

\eqref{E:B_q^n-2} and \eqref{E:B_q^n-4} imply that $\| V_{B_q^n}
\|_2^{1/n} \simeq \| V_{B_q^n} \|_4^{1/n}$. This fact, together with Lemma
\ref{T:Khintchine} and Proposition \ref{T:V_K,U_K}, imply that for any $1
\le p \le 4$ and $1 \le q \le \infty$,
\begin{eqnarray*}
\| V_{B_q^n} \|_p^{1/n} & \simeq & \sqrt{\alpha(1/q)} n^{-1/2}, \\
\| U_{B_q^n} \|_p^{1/n} & \simeq & \frac{1}{2} \sqrt{\alpha(1/q)}
    n^{-1/2}.
\end{eqnarray*}
$\alpha$ achieves its maximum over $[0,1]$ at $t = 0$, which implies that
for $1 \le p \le 4$ and $n$ large enough, $\| V_{B_q^n} \|_p$ and $\|
U_{B_q^n} \|_p$ are maximized over $1 \le q \le \infty$ when $q = \infty$.

Schmuckenschl\"ager \cite{Schmuck2} showed that $L_{B_q^n}$, and thus $\|
V_{B_q^n} \|_2$, is maximized by $q = \infty$ for all $n$. This suggests
the conjecture that $B_\infty^n$ is the body $K \in \K_s^n$ for which $\|
V_K \|_p$ and $\| U_K \|_p$ are maximized when $p$ is small enough. (Note
that the precise meaning of ``small enough" probably depends on the value
of $n$.) This may seem surprising initially, since a random crosspolytope
in $B_1^n$ can have full volume, whereas the maximum volume of a
crosspolytope in $B_\infty^n$ is much smaller when $n \ge 3$. The fact
that small moments are larger for $B_\infty^n$ reflects the fact that
there are many positions for a maximum volume crosspolytope in
$B_\infty^n$. On the other hand, since $\|V_{B_1^n}\|_\infty = 1$, it is
natural to conjecture that $\| V_K \|_p$ is maximized for some body which
is close to $B_1^n$ when $p$ is large enough.

We also note that for $1 \le q \le 2$, we have
$$\|V_{B_q^n}\|_\infty^{1/n}=\left(\frac{|B_1^n|}{|B_q^n|}\right)^{1/n}
    \simeq \frac{e^{1-1/q}}{q^{1/q} \Gamma(1 + \frac{1}{q})}
        n^{-\left(1 - \frac{1}{q}\right)}.$$
Therefore, $c_1 \le \sqrt{n} \|V_{B_2^n}\|_p^{1/n} \le c_2$ for $1 \le p
\le \infty$, where $c_1, c_2 > 0$ are absolute constants, but for $1 \le q
< 2$, the moment growth of $V_{B_q^n}$ is more complicated. This suggests
that in general, the asymptotic behavior of $\| V_K \|_p$ when $p$ and $n$
both increase without bound can depend strongly on the relationship
between $p$ and $n$.

\begin{prop} \label{T:asymp}
For any $0 < p < \infty$ and $K \in \K_s^n$, we have
$$\|V_K\|_p \le \left(1+\frac{p}{n}\right)^{-n/p}\|V_K\|_\infty.$$
\end{prop}
\begin{proof}
By integrating in spherical coordinates we obtain
\begin{eqnarray*}
\E V_K^p & = & \frac{1}{(n+p)^n |K|^n}
    \int\limits_{(S^{n-1})^n}
    \left(\prod_{i=1}^n \rho_K(\theta_i)\right)^n
    \frac{|\conv \{\pm \rho_K(\theta_1) \theta_1, \ldots,
    \pm \rho_K(\theta_n) \theta_n\}|^p}{|K|^p}
        \prod_{i=1}^n d\sigma(\theta_i) \\
& \le & \frac{\| V_K \|_\infty^p}{(n+p)^n |K|^n}
    \left( \int_{S^{n-1}} \rho_K(\theta)^n d\sigma(\theta) \right)^n
     =  \left(\frac{n}{n+p}\right)^n \| V_K \|_\infty^p,
\end{eqnarray*}
where $\rho_K(\theta) = \max \{ r > 0 : r \theta \in K \}$, and we have
used the fact that $\rho_K(\theta) \theta \in K$ for all $\theta \in
\R^n$.
\end{proof}

Note that $(1+1/t)^{-t} < 1$ for $t > 0$ and $\lim_{t \to \infty} (1 +
1/t)^{-t} = e^{-1} < 1$. Suppose we have a family $\{K_n \in \K_s^n : n
\in \N\}$ and let $p = p(n)$ so that $p = O(n)$. Then $\limsup_{n \to
\infty} (1 + p/n)^{-n/p} < 1$, so by Proposition \ref{T:asymp},
$$\limsup_{n \to \infty} \| V_{K_n} \|_p
    < \limsup_{n \to \infty} \| V_{K_n} \|_\infty.$$
Therefore, one can only have
$$\lim_{n \to \infty} \| V_{K_n} \|_p
    = \lim_{n \to \infty} \| V_{K_n} \|_\infty,$$
if $p$ grows faster than linearly with respect to $n$, for example, if $n
= o(p)$.

We now take up the proof of Proposition \ref{T:moment4}.

\begin{proof}[Proof of Proposition \ref{T:moment4}]
Let $S_n$ denote the group of permutations of $\{ 1, \ldots, n \}$. We
have
\begin{eqnarray}
\E V_K^{4}& = & \left( \frac{2^n}{(n! |K|)} \right)^{4} \frac{1}{|K|^n}
    \int_{K^n} [\det (x^1, \ldots, x^n)]^{4} dx^1\cdots dx^n \nonumber \\
& = & \frac{16^n}{(n!)^{4}|K|^{n + 4}} \int_{K^n}
    \prod_{i=1}^{4} \left(\sum_{ \tau_i \in S_n}\sgn (\tau_i)
    \prod_{j=1}^n x^j_{\tau_i(j)} \right) dx^1 \cdots dx^n \nonumber \\
& = & \frac{16^n}{(n!)^{4}|K|^{n + 4}}
    \sum_{\tau_1, \tau_2, \tau_3, \tau_{4} \in S_n}
    \sgn (\tau_1 \tau_2 \tau_3 \tau_{4}) \prod_{j=1}^n \int_K
    \prod_{i=1}^{4} x_{\tau_i(j)} dx \nonumber \\
& = & \frac{16^n}{(n!)^3 |K|^{n+4}}
    \sum_{\tau_1, \tau_2, \tau_3 \in S_n}
    \sgn ( \tau_1 \tau_2 \tau_3 )
    \prod_{j=1}^n \int_K x_j x_{\tau_1(j)} x_{\tau_2(j)} x_{\tau_3(j)} dx,
\label{E:moment4-1}
\end{eqnarray}
where in the last step we have used the fact that multiplication by
$\tau_4^{-1}$ permutes $S_n$.

We begin by defining $T_n$ to be the set of all triples $(\tau_1, \tau_2,
\tau_3) \in (S_n)^3$ such that for each $j = 1, \ldots, n$, (at least) one
of $\tau_1(j), \tau_2(j), \tau_3(j)$ is $j$ and the other two are equal.
Note first that if $K$ is symmetric with respect to reflections in the
coordinate hyperplanes, then the product in \eqref{E:moment4-1} is nonzero
only when $(\tau_1, \tau_2, \tau_3) \in T_n$. In this case,
\eqref{E:moment4-1} is simplified by the following fact.

\begin{claim} \label{T:simple}
If $(\tau_1, \tau_2, \tau_3) \in T_n$, then $\sgn(\tau_1 \tau_2 \tau_3) =
1$.
\end{claim}

Claim \ref{T:simple} can be proved by induction on $n$. With this,
\eqref{E:moment4-1} simplifies to
\begin{equation} \label{E:moment4-unc}
\E V_K^4 = \frac{16^n}{(n!)^3 |K|^{n+4}}
    \sum_{(\tau_1, \tau_2, \tau_3) \in T_n}
    \prod_{j=1}^n \int_K x_j x_{\tau_1(j)} x_{\tau_2(j)} x_{\tau_3(j)} dx
\end{equation}

Now if $K$ is also symmetric with respect to permutations of the
coordinates, then the integral expression in \eqref{E:moment4-unc} is
equal to $\int_K x_1^4 dx$ in the case that $\tau_1(j) = \tau_2(j) =
\tau_3(j) = j$, and equal to $\int_K x_1^2 x_2^2 dx$ otherwise. If we let
$d_{n,k}$ denote the number of triples $(\tau_1, \tau_2, \tau_3) \in T_n$
such that $\tau_1(j) = \tau_2(j) = \tau_3(j) = j$ for exactly $k$ values
of $j \in \{ 1, \ldots, n \}$, then we have
\begin{eqnarray*}
\E V_K^4 & = & \frac{16^n}{(n!)^3 |K|^{n+4}}
    \sum_{k = 0}^n d_{n,k} \left(\int_K x_1^4 dx \right)^k
    \left(\int_K x_1^2 x_2^2 dx \right)^{n-k} \\
    & = & \frac{(16 \int_K x_1^2 x_2^2 dx)^n}{(n!)^3 |K|^{n+4}}
    \sum_{k = 0}^n d_{n,k} \left( \frac{\int_K x_1^4 dx}
        {\int_K x_1^2 x_2^2 dx} \right)^k.
\end{eqnarray*}

Proposition \ref{T:moment4} will now follow from the following.

\begin{claim} \label{T:app}
For any $n \in \N$ and $t \in \R$, we have
$$\sum_{k = 0}^n d_{n,k} t^k = \frac{(n+2)!}{2} \varphi(t-3,n),$$
where $\varphi(t,n)$ is as defined in \eqref{E:phi}.
\end{claim}

To prove this, we first note that if we define $d_n = d_{n,0}$ for $n \in
\N$ and $d_0 = 1$, then $d_{n,k} = \binom{n}{k} d_{n-k}$. The sequence $\{
d_n : n \in \N \cup \{0\} \}$ satisfies the following recurrence relation.

\begin{claim} \label{T:recurrence}
For each $n \in \N$, $d_{n+1} = n (d_n + 3 d_{n-1})$.
\end{claim}
\begin{proof}
Let $\tilde{T}_n = \{ (\tau_1, \tau_2, \tau_3) \in T_n : \tau_1(j),
\tau_2(j), \tau_3(j)$ are never all equal$\}$. Then $d_n = |\tilde{T}_n|$
($| \cdot |$ denotes the cardinality of a finite set) . Let $(\tau_1,
\tau_2, \tau_3) \in \tilde{T}_{n+1}$, and let $\omega \in S_{n+1}$ be as
defined by the condition that for each $j = 1, \ldots, n+1$, we have $\{
\tau_1(j), \tau_2(j), \tau_3(j) \} = \{ j, \omega(j), \omega(j) \}$ (as
multisets). Note $\omega$ has no fixed points. Now $\tilde{T}_n$ can be
partitioned according to the $n$ possible values of $\omega(n)$; the
number of triples in $\tilde{T}_n$ corresponding to each of these values
of $\omega(n)$ is equal. Therefore
$$d_{n+1} = n \bigl|\{ (\tau_1, \tau_2, \tau_3) \in \tilde{T}_{n+1} :
        \omega(n) = n+1 \}\bigr|.$$
The set appearing in this expression can be further partitioned according
to the value of $\omega(n+1)$. If $\omega(n+1) = n$, then one of $\tau_1,
\tau_2, \tau_3$ fixes both $n$ and $(n+1)$ and the other two transpose
them; therefore each of $\tau_1, \tau_2, \tau_3$ restricts to a
permutation of $\{ 1, \ldots, n-1 \}$ in such a way that the restrictions
form a triple in $\tilde{T}_{n-1}$. Furthermore, each triple
$(\tilde{\tau}_1, \tilde{\tau}_2, \tilde{\tau}_3) \in \tilde{T}_{n-1}$
results in this way from exactly three triples $(\tau_1, \tau_2, \tau_3)
\in T_{n+1}$ such that $\omega(n) = n+1$ and $\omega(n+1) = n$. Therefore
there are $3 d_{n-1}$ triples such that $\omega(n) = n+1$ and $\omega(n+1)
= n$.

On the other hand, if $\omega(n) = n+1$ and $\omega(n+1) \ne n$, we can
define $(\tilde{\tau}_1, \tilde{\tau}_2, \tilde{\tau_3}) \in \tilde{T}_n$
by
$$ \tilde{\tau_i}(j) = \begin{cases}
     n        & \text{ if } j = n \text{ and } \tau_i(n+1) = n+1,   \\
    \tau_i(n) & \text{ if } j = n \text{ and } \tau_i(n+1) \ne n+1, \\
    \tau_i(j) & \text{ if } j < n.
    \end{cases}$$
One can easily verify that this defines a bijection between $\{ (\tau_1,
\tau_2, \tau_3) \in \tilde{T}_{n+1} : \omega(n) = n+1$ and $\omega(n+1)
\ne n \}$ and $\tilde{T}_n$. Therefore there are $d_n$ triples $(\tau_1,
\tau_2, \tau_3) \in \tilde{T}_{n+1}$ such that $\omega(n) = n+1$ and
$\omega(n+1) \ne n$.
\end{proof}

\begin{proof}[Proof of Claim \ref{T:app}]
Let $g$ be the exponential generating function of the sequence $\{ d_n : n
\in \N \cup \{0\} \}$, defined by
$$g(t) = \sum_{n = 0}^\infty d_n \frac{t^n}{n!}.$$
The recurrence in Claim \ref{T:recurrence} and the initial conditions $d_0
= 1$, $d_1 = 0$, imply that $g$ satisfies the differential equation
$$ g'(t) = \frac{3t}{1 - t} g(t)$$
with initial condition $g(0) = 1$, which has the solution
$$g(t) = \frac{e^{-3t}}{(1 - t)^3}.$$
Now we have
\begin{eqnarray*}
\sum_{n = 0}^\infty \left( \sum_{k = 0}^n d_{n,k} t^k \right)
    \frac{u^n}{n!}
& = & \sum_{n = 0}^\infty \left( \sum_{k = 0}^n \binom{n}{k} d_{n-k} t^k
    \right) \frac{u^n}{n!}
    =  \left( \sum_{n = 0}^\infty d_n \frac{u^n}{n!} \right)
    \left( \sum_{n = 0}^\infty t^n \frac{u^n}{n!} \right) \\
& = & \frac{e^{-3u}}{(1 - u)^3} e^{tu} = \frac{1}{(1 - u)^3} e^{(t-3)u} \\
& = & \left( \frac{1}{2} \sum_{n = 0}^\infty (n+1)(n+2) u^n \right)
    \left( \sum_{n = 0}^\infty \frac{(t - 3)^n}{n!} u^n \right) \\
& = & \frac{n!}{2} \sum_{n = 0}^\infty \left( \sum_{k = 0}^n
    \frac{(n-k+1)(n-k+2)}{k!} (t-3)^k \right) \frac{u^n}{n!},
\end{eqnarray*}
so that
$$\sum_{k = 0}^n d_{n,k} t^k = \frac{n!}{2} \sum_{k = 0}^n
    \frac{(n+1-k)(n+2-k)}{k!} (t-3)^k.$$
The remainder of the proof is elementary calculation.
\end{proof}

The proof of Claim \ref{T:app} follows the outline of a similar
calculation, with the numbers $d_n$ replaced by the derangement numbers
$D_n$, shown to the author by A.~de Acosta \cite{Acosta}. This completes
the proof of Proposition \ref{T:moment4}.
\end{proof}

\section{Random crosspolytopes in $B_2^n$} \label{S:euclidean}
We now specialize to the case $K = B_2^n$, in which it is possible to
compute explicitly all of the moments of $V_K$. In fact, with little more
effort it is possible to derive a more general result.

\begin{prop}\label{T:symm-miles}
Let $1\le k\le n$ and $0\le s \le k$. Consider $k$ independent random
points in $\R^n$, $s$ of which are uniformly distributed in $S^{n-1}$ and
$k-s$ of which are uniformly distributed in $B_2^n$. Let $V=V_{n,k,s}$
denote the $k$-dimensional volume of the symmetric convex hull of the
random points. Then for $p \ge 0$,
$$\E V^p = \left(\frac{2^k}{k!}\right)^p \left(1+\frac{p}{n}\right)^s
    \left(\frac{\Gamma(1+\frac{n}{2})}{\Gamma(1+\frac{n+p}{2})}\right)^k
    \prod_{i=1}^k\frac{\Gamma(\frac{n-k+i+p}{2})}{\Gamma(\frac{n-k+i}{2})}.
$$
\end{prop}

Proposition \ref{T:symm-miles} is the symmetric analogue of a similar
result due to Miles \cite{Miles}. The most interesting special case is
when $k=n$ and $s=0$, in which we have $V_{n,n,0} = |B_2^n| V_{B_2^n}$. By
the above formula, we have
$$\E V_{B_2^n}^p
    = \left(\frac{2^n \Gamma(1+\frac{n}{2})}{n!\pi^{n/2}}\right)^p
    \left(\frac{\Gamma(1+\frac{n}{2})}{\Gamma(1+\frac{n+p}{2})}\right)^n
    \prod_{i=1}^n\frac{\Gamma(\frac{i+p}{2})}{\Gamma(\frac{i}{2})}.$$

\begin{proof}[Proof of Proposition \ref{T:symm-miles}:] If
$x_1,\ldots,x_k\in\R^n$ and $A=A(x_1,\ldots,x_k)$ is the $n\times k$
matrix with columns $x_1,\ldots,x_k$, then we have
$$ v_k(x_1,\ldots,x_k) = \vol_k(\conv\{\pm x_1,\ldots,\pm x_k\})
    = \frac{2^k}{k!}|\det A^t A|^{1/2}.$$
Using the homogeneity of this quantity as a function of each point, we
proceed by first integrating out the radial dependence of the points which
are uniformly distributed in the ball, then transforming the resulting
integrals over spheres in the usual way to Gaussian integrals. In this way
we obtain
\begin{eqnarray*}
\E V^p & = & \frac{1}{|B_2^n|^{k-s} \sigma(S^{n-1})^s}
    \int\limits_{(S^{n-1})^s} \int\limits_{(B_2^n)^{k-s}}
    v_k(\omega_1,\ldots,\omega_s,x_{s+1},\ldots,x_k)^p
    d\sigma(\omega_1)\cdots d\sigma(\omega_s) dx_{s+1}\cdots dx_k \\
& = & \left(\frac{2^{k/2}}{k!}\right)^p \left(1+\frac{p}{n}\right)^s
    \left(\frac{\Gamma(1+\frac{n}{2})}{\Gamma(1+\frac{n+p}{2})}\right)^k
    \E |\det G^t G|^{p/2},
\end{eqnarray*}
where $G = G_{n, k}$ is an $n \times k$ random matrix whose entries are
independent $N(0,1)$ random variables. The proposition now follows from
the classical identity
$$\E |\det G^t G|^p = 2^{kp}\prod_{i=1}^k
    \frac{\Gamma(p+\frac{n-k+i}{2})}{\Gamma(\frac{n-k+i}{2})}.$$
A simple proof of this identity may be given, e.g., by using the
representation of $G$ in \cite{Silver}.
\end{proof}


\section*{Acknowledgments}
This paper is part of the author's Ph.D.~thesis, written under the
supervision of Profs.~S. Szarek and E.~Werner. The author wishes to thank
Profs.~E.~Werner and C.~Sch\"utt for many valuable discussions.


\end{document}